\documentstyle[11pt,amstex,amssymb]{amsart}
\newtheorem{thm}{Theorem}[section]

\newtheorem{lemma}[thm]{Lemma}

\newcommand{\proof}{\noindent{\it Proof.}\enspace}

\def\bT{{\Bbb T}}
\def\Ric{{\rm Ric}}
\def\cM{{\cal M}}
\def\bR{{\Bbb R}}
\def\bbbr{{\Bbb R}}
\def\Hess{{\rm Hess}}
\def\Tr{{\rm Tr}}
\def\SU{{\rm SU}}
\def\bN{{\Bbb N}}
\def\im{{\rm i}}
\begin{document}

\title[A free TCI on the circle: a new proof]
{A free analogue of the transportation cost inequality on the circle}

\author[F. Hiai]{Fumio Hiai$\,^{1,2}$}
\address{Graduate School of Information Sciences,
Tohoku University, Aoba-ku, Sendai 980-8579, Japan}
\author[D. Petz]{D\'enes Petz$\,^{1,3}$}
\address{Alfr\'ed R\'enyi Institute of Mathematics, Hungarian 
Academy of Sciences, H-1053 Budapest, Re\'altanoda u. 13-15, Hungary}

\thanks{$^1\,$Supported in part by Japan-Hungary JSPS-HAS Joint Project.}

\thanks{$^2\,$Supported in part by Strategic Information and Communications
R\&D Promotion Scheme of MPHPT}

\thanks{$^3\,$Supported in part by OTKA T032662.}


\maketitle
\centerline{\it Dedicated to Professor Andr\'as Pr\'ekopa on
the occasion of his 75th birthday}
\bigskip
\bigskip
\bigskip

The relative entropy
$$
S(\mu,\nu):=\cases
\displaystyle\int\log{d\mu\over d\nu}\,d\mu & \text{if $\mu\ll\nu$}, \\
\ +\infty & \text{otherwise}.
\endcases
$$
and the Wasserstein distance are useful distances between
measures. For probability measures $\mu$ and $\nu$ on the Euclidean space
$\bR^n$, the latter is defined as
$$
W(\mu,\nu):=\inf_{\pi\in\Pi(\mu,\nu)}
\sqrt{\iint {1\over2}d(x,y)^2\,d\pi(x,y)},
$$
where $d(x,y)=\|x-y\|_2$ and $\Pi(\mu,\nu)$ denotes the set of all 
probability measures on $\bbbr^n\times\bbbr^n$ with marginals $\mu$ and 
$\nu$, i.e., $\pi(\,\cdot\times\bbbr^n)=\mu$ and $\pi(\bbbr^n\times\cdot\,)
=\nu$. 

The transportation cost inequality (TCI) obtained by M.~Talagrand \cite{Ta}
is 
$$
W(\mu,\nu)\le\sqrt{S(\mu,\nu)}\, ,
$$ 
where $\nu$ is the standard Gaussian measure and $\mu$ is any probability 
measure on $\bR^n$. Recently Talagrand's inequality and its counterpart,
the logarithmic Sobolev inequality (LSI) have received a lot of attention
and they have been extended from the Euclidean spaces to Riemannian manifolds.
(Contrary to Talagrand's inequality, the LSI gives an upper bound for the
relative entropy.) It was shown by F.~Otto and  C.~Villani \cite{OV} that in
the Riemannian manifold setting the TCI  follows from the LSI due to
D.~Bakry and M.~Emery \cite{BE}.

On the other hand, Ph.~Biane and D.~Voiculescu \cite{BV} proved the free
analogue of Talagrand's TCI for compactly supported measures on the real line.
They replaced the relative entropy with the relative free entropy and the 
Gaussian measure with the semicircular law. Based on the method of 
random matrix approximation, Biane \cite{Bi} proved the free LSI for 
measures on the real line, and we made a slight generalization of Biane
and Voiculescu's free TCI \cite{HPU1}. We also obtained the free TCI and 
LSI for measures on the unit circle using large deviation results for
special unitary matrices and the differential geometry of $\SU(n)$ 
\cite{HPU1, HPU2}.

Recently M.~Ledoux \cite{Le} used the random matrix method to obtain the
free analogue  of the Pr\'ekopa-Leindler inequality on the real line.
From this together with the Hamilton-Jacobi approach, he also gave
different proofs of the free LSI and TCI for measures on $\bR$. The aim of
the present notes is to give a new proof of the free TCI for measures on
the circle following Ledoux's idea. In this way we do not need the large
deviation technique but we establish a kind of free analogue of 
Pr\'ekopa-Leindler inequality on the circle.

\section{The Pr\'ekopa-Leindler inequality on a Riemannian manifold}

Let $M$ be a complete, connected, $n$-dimensional Riemannian manifold with
the volume measure $dx$ and the geodesic distance $d(x,y)$ for $x,y\in M$.
For $0<\theta<1$ define
$$
Z_\theta(x,y):=\Bigl\{z\in M:
d(x,z)=\theta d(x,y),\,d(z,y)=(1-\theta)d(x,y)\Bigr\},
$$
which is the locus of points playing the role of $(1-\theta)x+\theta y$.
In this section we first present a result of Cordero-Erausquin, McCann and 
Schmuckenschl\"ager, which is an extension of the {Pr\'ekopa-Leindler 
inequality} to the Riemannian manifold setting. Then we show
that this results implies the TCI on a Riemannian manifold under some
conditions (slightly stronger than the Bakry-Emery criterion).

\begin{thm}\label{T-1.1} {\rm (\cite[Corollary 1.2]{CMS})}\quad
Assume that $\Ric(M)\ge(n-1)k$ holds for some $k\in\bR$ where $\Ric(M)$ is
the Ricci curvature of $M$. Let $f,g,h:M\to[0,\infty)$ be Borel measurable
functions and fix $0<\theta<1$. Assume that
$$
h(z)\ge\biggl({S_k(d)\over S_k((1-\theta)d)^{1-\theta}
S_k(\theta d)^\theta}\biggr)^{n-1}f(x)^{1-\theta}g(y)^\theta
$$
holds for every $x,y\in M$, $z\in Z_\theta(x,y)$ and $d:=d(x,y)$, where
$$
S_k(d):=\cases
\sin(\sqrt{k}d)/\sqrt{k}d & \text{if $k>0$}, \\
\ 1 & \text{if $k=0$}, \\
\sinh(\sqrt{-k}d)/\sqrt{-k}d & \text{if $k<0$}.
\endcases
$$
Then
$$
\int_Mh(x)\,dx\ge\biggl(\int_Mf(x)\,dx\biggr)^\theta
\biggl(\int_Mg(x)\,dx\biggr)^{1-\theta}.
$$
\end{thm}

Here it is worth noting a known result: If $\Ric(M)\ge(n-1)k$ with $k>0$, 
then the diameter of $M$ is at most $\pi/\sqrt{k}$ (see \cite[1.26]{CE}). 

Write
$$
\Phi_\theta(d):=(n-1)\bigl(\log S_k(d)-(1-\theta)\log S_k((1-\theta)d)
-\theta\log S_k(\theta d)\bigr).
$$
Let $\cM(M)$ denote the set of probability Borel measures on $M$. Let
$\nu\in\cM(M)$ be given by $d\nu:={1\over Z}e^{-Q(x)}\,dx$ with a Borel
function $Q:M\to\bR$ and a normalization constant $Z$. Write
$$
R_\theta(z;x,y):=Q(z)-(1-\theta)Q(x)-\theta Q(y).
$$
Then, the above theorem is rephrased as follows: 

If $u,v,w:M\to\bR$ are Borel functions and
$$
w(z)\ge(1-\theta)u(x)+\theta v(y)+R_\theta(z;x,y)+\Phi_\theta(d)
$$
holds for every $x,y\in M$, $z\in Z_\theta(x,y)$ and $d:=d(x,y)$, then
\begin{equation}\label{F-1.1}
\log\int_Me^{w(x)}\,d\nu(x)\ge(1-\theta)\log\int_Me^{u(x)}\,d\nu(x)
+\theta\log\int_Me^{v(x)}\,d\nu(x).
\end{equation}

The following transportation cost inequality in the Riemannian setting was
shown in \cite{OV} based on \cite{BE}.

\begin{thm}\label{T-1.2} {\rm (\cite{BE} and \cite{OV})}\quad
Let $\nu\in\cM(M)$ be given by $d\nu(x):={1\over Z}e^{-Q(x)}\,dx$ with a
$C^2$ function $Q:M\to\bR$. If the Bakry and Emery criterion
\begin{equation}\label{F-1.2}
\Ric(M)+\Hess(Q)\ge\rho I_n
\end{equation}
is satisfied with a constant $\rho>0$, then
$$
W(\mu,\nu)\le\sqrt{{1\over\rho}S(\mu,\nu)},
\qquad\mu\in\cM(M).
$$
\end{thm}

Now, we assume the following condition slightly stronger than
\eqref{F-1.2}:
\begin{equation}\label{F-1.3}
\Ric(M)\ge\alpha I_n\quad{\rm and}\quad\Hess(Q)\ge\beta I_n
\end{equation}
for some constants $\alpha>0$, $\beta\in\bR$ with $\alpha+\beta=\rho>0$.
Our goal in this section is to prove that Theorem \ref{T-1.1} implies
Theorem \ref{T-1.2} under the assumption \eqref{F-1.3}. 

We use the celebrated variational formula (or the {Monge-Kantorovich 
duality}) for the Wasserstein distance (see \cite{Vi}):
\begin{eqnarray}\label{F-1.4}
&&\rho W(\mu,\nu)^2=\sup\biggl\{\int_Mf(x)\,d\mu(x)-\int_Mg(x)\,d\nu(x):
\nonumber\\
&&\hskip2.5cm f,g\in C_b(M),\,f(x)\le g(y)+{\rho\over2}d(x,y)^2,
\,x,y\in M\biggr\}
\end{eqnarray}
where $\rho>0$. The variational expression for the relative entropy is
also useful:
\begin{equation}\label{F-1.5}
S(\mu,\nu)=\sup\biggl\{\int_Mf(x)\,d\mu(x)-\log\int_Me^{f(x)}\,d\nu(x)
:f\in C_b(M)\biggr\}.
\end{equation}
Furthermore, we need the Taylor expansion of $\log S_k(d)$:

\begin{lemma}\label{L-1.3}
$$
\log S_k(d)=-\sum_{j=1}^\infty c_j(kd^2)^j
$$
with $c_1={1\over6}$ and $c_j>0$ for all $j\ge1$.
\end{lemma}

\proof
Set $f(x):=\log{\sin x\over x}$; then we have
$$
f'(x)=\cot x-{1\over x}=\sum_{m=1}^\infty{2x\over x^2-(m\pi)^2}
$$
because of a well-known expansion of $\cot x$. Since
$$
f(x)=-\sum_{m=1}^\infty{1\over m\pi}
\biggl({1\over1-{x\over m\pi}}-{1\over1+{x\over m\pi}}\biggr)
=-\sum_{m=1}^\infty\sum_{j=1}^\infty{2\over(m\pi)^{2j}}x^{2j-1},
$$
we get
$$
f^{(2j)}(0)=-{2(2j-1)!\over\pi^{2j}}\sum_{m=1}^\infty{1\over m^{2j}}
$$
so that
$$
c_j=-{f^{(2j)}(0)\over(2j)!}
={1\over j\pi^{2j}}\sum_{m=1}^\infty{1\over m^{2j}}>0.
$$
\qed

\bigskip
Since $\Ric(M)\ge(n-1)k$ with $k={\alpha\over n-1}$ due to the assumption
$\Ric(M)\ge\alpha I_n$ in \eqref{F-1.3}, we get by Lemma \ref{L-1.3}
\begin{eqnarray}\label{F-1.6}
\Phi_\theta(d)&=&-(n-1)\sum_{j=1}^\infty
\bigl(1-(1-\theta)^{2j+1}-\theta^{2j+1}\bigr)c_j(kd^2)^j \nonumber \\
&\le&-(n-1)\bigl(1-\theta^3-(1-\theta)^3\bigr){\alpha\over6(n-1)}d^2
\nonumber\\
&=&-{\alpha\theta(1-\theta)\over2}d^2.
\end{eqnarray}
For each $x,y\in M$ let $z(t)$ ($0\le t\le1$) be a geodesic curve joining
$x,y$ with $d(x,z(t))=td(x,y)$. Since the assumption
$\Hess(Q)\ge\beta I_n$ in \eqref{F-1.3} gives
$$
{d^2\over dt^2}Q(z(t))\ge\beta d(x,y)^2,\qquad 0\le t\le1,
$$
we get
\begin{eqnarray}\label{F-1.7}
R_\theta(z(\theta);x,y)
&=&Q(z(\theta))-\theta Q(z(0))-(1-\theta)Q(z(1)) \nonumber\\
&\le&-{\beta\theta(1-\theta)\over2}d(x,y)^2.
\end{eqnarray}
Hence, by \eqref{F-1.6} and \eqref{F-1.7} we have
$$
R_\theta(z;x,y)+\Phi_\theta(d(x,y))
\le-{\rho\theta(1-\theta)\over2}d(x,y)^2
$$
for every $x,y\in M$ and $z\in Z_\theta(x,y)$.

Now, let $f,g\in C_b(M)$ be such that
$$
f(x)\le g(y)+{\rho\over2}d(x,y)^2,\qquad x,y\in M.
$$
Set $u:=\theta f$, $v:=-(1-\theta)g$ and $w:=0$. Then
\begin{eqnarray*}
&&(1-\theta)u(x)+\theta v(y)+R_\theta(z;x,y)+\Phi_\theta(d(x,y)) \\
&&\qquad\le\theta(1-\theta)\Bigl\{f(x)-g(y)-{\rho\over2}d(x,y)^2\Bigr\}
\le0=w(z)
\end{eqnarray*}
for every $x,y\in M$ and $z\in Z_\theta(x,y)$. Hence Theorem \ref{T-1.1}
(the rephrased version \eqref{F-1.1}) yields
$$
\log\int_Me^{\theta f(x)}\,d\nu(x)
+{\theta\over1-\theta}\log\int_Me^{-(1-\theta)g(x)}\,d\nu(x)\le0.
$$
Letting $\theta\nearrow1$ gives
$$
\log\int_Me^{f(x)}\,d\nu(x)-\int_Mg(x)\,d\nu(x)\le0
$$
so that
$$
\int_Mf(x)\,d\mu(x)-\int_Mg(x)\,d\nu(x)
\le\int_Mf(x)\,d\mu(x)-\log\int_Me^{f(x)}\,d\nu(x)\le S(\mu,\nu)
$$
thanks to \eqref{F-1.5}. Finally, we apply \eqref{F-1.4} to obtain
$$
\rho W(\mu,\nu)^2\le S(\mu,\nu).
$$
\qed

\section{Free TCI on the circle}
\setcounter{equation}{0}

Let $Q:\bT\to\bR$ be a continuous function. The {weighted energy
integral} associated with $Q$ is defined by
$$
E_Q(\mu):=-\Sigma(\mu)+\int_\bT Q(\zeta)\,d\mu(\zeta)
\quad\mbox{for $\mu\in\cM(\bT)$},
$$
which admits a unique minimizer $\nu_Q\in\cM(\bT)$ (see \cite{ST}). Set
$B(Q):=-E_Q(\nu_Q)$ and define the {relative free entropy} with respect
to $Q$ by
$$
\widetilde\Sigma_Q(\mu):=-\Sigma(\mu)+\int_\bT Q(\zeta)\,d\mu(\zeta)+B(Q)
\quad\mbox{for $\mu\in\cM(\bT)$}.
$$
It is known (\cite[Theorem 2.1]{HPU2}, also \cite[Chap.\ 5]{HP}) that
$\widetilde\Sigma_Q(\mu)$ is the rate function of the large deviation
principle (in the scale
$1/N^2$) for the empirical eigenvalue distribution of the special unitary
random matrix
$$
d\lambda_N^{\rm SU}(Q)(U):={1\over Z_N^{\rm SU}(Q)}
\exp\bigl(-N\Tr_N(Q(U))\bigr)\,dU,
$$
where $dU$ is the Haar probability measure on the special unitary group
$\SU(N)$ of order $N$, $Q(U)$ for $U\in\SU(N)$ is defined via functional
calculus and $\Tr_N$ is the usual trace on the $N\times N$ matrices.

The Wasserstein distance $W(\mu,\nu)$ between $\mu,\nu\in\cM(\bT)$ is
defined with respect to the angular distance (i.e., the geodesic
distance). The following is the free TCI for measures on $\bT$ proven in
\cite{HPU1}. The aim of this section is to re-prove this by using the
method of Ledoux \cite{Le}.

\begin{thm}\label{T-2.1} {\rm (\cite[Theorem 2.7]{HPU1})}\quad
Let $Q:\bT\to\bR$ be a continuous function. If there exists a
constant $\rho>-{1\over2}$ such that $Q(e^{\im t})-{\rho\over2}t^2$
is convex on $\bR$, then
$$
W(\mu,\nu_Q)\le\sqrt{{2\over1+2\rho}\widetilde\Sigma_Q(\mu)},
\qquad\mu\in\cM(\bT).
$$
\end{thm}

We introduce the {relative free pressure} with respect to $Q$ by
$$
j_Q(f):=\sup\biggl\{\int_\bT f\,d\mu-\widetilde\Sigma_Q(\mu)
:\mu\in\cM(\bT)\biggr\}\quad\mbox{for $f\in C_\bR(\bT)$.}
$$
It is known (\cite{Le} and \cite{HMP}) that
\begin{eqnarray}\label{F-2.1}
j_Q(f)&=&E_Q(\nu_A)-E_{Q-f}(\nu_{Q-f}) \nonumber\\
&=&\lim_{N\to\infty}{1\over N^2}\log\int_{\SU(N)}
\exp\bigl(N\Tr_N(f(U))\bigr)\,d\lambda_N^\SU(Q)(U).
\end{eqnarray}
For $N\in\bN$ and $U,V,W\in\SU(N)$ write
$$
R_{\theta,N}(W;U,V):=\Tr_N(Q(W))-(1-\theta)\Tr_N(Q(U))-\theta\Tr_N(Q(V)).
$$
The next lemma is a sort of free analogue of Pr\'ekopa-Leindler-Ledoux
inequality on the circle.

\begin{lemma}\label{L-2.2}
Let $f,g,h:\bT\to\bR$ be Borel functions and fix $0<\theta<1$. Assume that
\begin{eqnarray*}
\Tr_N(h(W))&\ge&(1-\theta)\Tr_N(f(U))+\theta\Tr_N(g(V)) \\
&&\qquad\qquad+R_{\theta,N}(W;U,V)-{\theta(1-\theta)\over4}d(U,V)^2
\end{eqnarray*}
holds for every $N\in\bN$, $U,V\in\SU(N)$ and $W\in Z_\theta(U,V)$. Then
\begin{equation}\label{F-2.2}
j_Q(h)\ge(1-\theta)j_Q(f)+\theta j_Q(g).
\end{equation}
\end{lemma}

\proof
Since $\dim(\SU(N))=N^2-1$ and $\Ric(\SU(N))={N\over2}$, $\Phi_\theta$
defined for $M=\SU(N)$ satisfies
$$
\Phi_\theta(d)\le-{N\theta(1-\theta)\over4}d^2
$$
thanks to \eqref{F-1.6}. Hence, for each $N\in\bN$, the assumption of the
lemma gives
\begin{eqnarray*}
N\Tr_N(h(W))&\ge&(1-\theta)N\Tr_N(f(U))+\theta N\Tr_N(g(V)) \\
&&\qquad\qquad+NR_{\theta,N}(W;U,V)+\Phi_\theta(d(U,V))
\end{eqnarray*}
for every $U,V\in\SU(N)$ and $W\in Z_\theta(U,V)$. Theorem \ref{T-1.1}
(the rephrased version \eqref{F-1.1}) can be applied to
$\nu:=\lambda_N^\SU(Q)$, $u:=N\Tr_N(f(\cdot))$, $v:=N\Tr_N(g(\cdot))$ and
$w:=N\Tr_N(h(\cdot))$; hence we have
\begin{eqnarray*}
&&\log\int_{\SU(N)}\exp\bigl(N\Tr_N(h(U))\bigr)\,d\lambda_N^\SU(Q)(U) \\
&&\qquad\ge(1-\theta)
\log\int_{\SU(N)}\exp\bigl(N\Tr_N(f(U))\bigr)\,d\lambda_N^\SU(Q)(U) \\
&&\qquad\qquad+\theta
\log\int_{\SU(N)}\exp\bigl(N\Tr_N(g(U))\bigr)\,d\lambda_N^\SU(Q)(U),
\end{eqnarray*}
implying the inequality \eqref{F-2.2} thanks to \eqref{F-2.1}.\qed

\bigskip
The assumption of the lemma is apparently too much; so the above must not be
the optimal form of the free Brunn-Minkowski inequality on $\bT$.
Nevertheless, it is enough to prove Theorem \ref{T-2.1}.

For each
$N\in\bN$ and $U\in\SU(N)$ set $\Psi(U):=\Tr_N(Q(U))$. Using a certain
regularization technique as in \cite{HPU1}, we may assume that $Q$ is a
harmonic function in a neighborhood of the unit disk. Then, it was shown
in \cite[Lemma 1.3\,(ii)]{HPU1} that the convexity assumption of $Q$
implies $\Hess(\Psi)\ge\rho I_{N^2-1}$. This gives as in \eqref{F-1.7}
\begin{equation}\label{F-2.3}
R_{\theta,N}(W;U,V)\le-{\rho\theta(1-\theta)\over2}d(U,V)^2
\end{equation}
for every $U,V\in\SU(N)$ and $W\in Z_\theta(U,V)$. Now, let
$f,g\in C(\bT)$ be such that
\begin{equation}\label{F-2.4}
f(\zeta)\le g(\eta)+{1+2\rho\over4}d(\zeta,\eta)^2,
\qquad\zeta,\eta\in\bT.
\end{equation}
Define the optimal matching distance on $\bT^N$ by
$$
\delta(\zeta,\eta):=\min_{\sigma\in S_N}
\sqrt{\sum_{i=1}^Nd(\zeta_i,\eta_{\sigma(i)})^2}
$$
for $\zeta=(\zeta_1,\dots,\zeta_N),\eta=(\eta_1,\dots,\eta_N)\in\bT^N$.
For $U\in\SU(N)$ let $\lambda(U):=(\lambda_1(U),\dots,\lambda_N(U))$
denote the element of $\bT^N$ consisting of the eigenvalues of $U$ with
multiplicities and in counter-clockwise order. It immediately follows from
\eqref{F-2.4} that
$$
\Tr_N(f(U))\le\Tr_N(g(V))+{1+2\rho\over4}\delta(\lambda(U),\lambda(V))^2,
\qquad U,V\in\SU(N).
$$
Since $\delta(\lambda(U),\lambda(V))\le d(U,V)$ as shown in
\cite[(2.11)]{HPU1}), this gives
\begin{equation}\label{F-2.5}
\Tr_N(f(U))\le\Tr_N(g(V))+{1+2\rho\over4}d(U,V)^2,
\qquad U,V\in\SU(N).
\end{equation}
Set $\tilde f:=\theta f$, $\tilde g:=-(1-\theta)g$ and $\tilde h:=0$.
Then, for $U,V\in\SU(N)$ and $W\in Z_\theta(U,V)$, by \eqref{F-2.3} and
\eqref{F-2.5} we get
\begin{eqnarray*}
&&(1-\theta)\Tr_N(\tilde f(U))+\theta\Tr_N(\tilde g(V))
+R_{\theta,N}(W;U,V)-{\theta(1-\theta)\over4}d(U,V)^2 \\
&&\qquad\le\theta(1-\theta)\biggl(\Tr_N(f(U))-\Tr_N(g(V))
-{1+2\rho\over4}d(U,V)^2\biggr) \\
&&\qquad\le0=\Tr_N(\tilde h(W)).
\end{eqnarray*}
Hence, the assumption of Lemma \ref{L-2.2} is satisfied so that we have
$$
(1-\theta)j_Q(\theta f)+\theta j_Q(-(1-\theta)g)\le j_Q(0)=0.
$$
For every $\mu\in\cM(\bT)$, by definition of $j_Q$, this implies
$$
(1-\theta)\biggl(\int_\bT\theta f\,d\mu-\widetilde\Sigma_Q(\mu)\biggr)
+\theta\biggl(\int_\bT(-(1-\theta)g)\,d\nu_Q
-\widetilde\Sigma_Q(\nu_Q)\biggr)\le0
$$
so that, thanks to $\widetilde\Sigma_Q(\nu_Q)=0$,
$$
\theta\biggl(\int_\bT f\,d\mu-\int_\bT f\,d\nu_Q\biggr)
\le\widetilde\Sigma_Q(\mu).
$$
Letting $\theta\nearrow1$ gives
$$
\int_\bT f\,d\mu-\int_\bT f\,d\nu_Q\le\widetilde\Sigma_Q(\mu).
$$
Using \eqref{F-1.4} we obtain
$$
{1+2\rho\over2}W(\mu,\nu_Q)^2\le\widetilde\Sigma_Q(\mu).
$$
\qed

\bigskip
It turns out that the bound ${2/(1+2\rho)}$ of our free TCI on $\bT$ 
cannot be improved even if we use the Riemannian
Pr\'ekopa-Leindler inequality from \cite{CMS}. This suggests the best
possibility of the bound.

\bigskip\noindent
{\bf Acknowledgments.}\enspace
We are grateful to Professor M.~Ledoux for sending us his preprint
\cite{Le}, to Professor Y.~Ueda for suggesting us the proof of free TCI by
using \cite{CMS}, and to Dr.\ A.~Andai for the proof of Lemma \ref{L-1.3}.

\end{document}